\documentclass[12pt]{article}
\usepackage{amsmath,amssymb,amsthm,times,authblk,mathtools,url,cite}
\usepackage[colorlinks=true, linkcolor=blue, urlcolor=black, citecolor = blue, anchorcolor = blue, linktocpage]{hyperref}
\usepackage{color}
\usepackage{orcidlink}
\usepackage[normalem]{ulem}
\usepackage{mathrsfs}
\newcommand\blue{\textcolor{blue}}

\newcommand{\comment}[1]{}

\newtheorem{theorem}{Theorem}
\newtheorem{lemma}{Lemma}
\newtheorem{corollary}{Corollary}

\theoremstyle{remark}
\newtheorem{remark}{Remark}

\allowdisplaybreaks 

\title{Quasi-orthogonal $q$-Meixner polynomials}
\author{Kerstin Jordaan$^{a,b}\ \orcidlink{0000-0002-1675-5366},$ Pinaki Prasad Kar$^{c}\ \orcidlink{0000-0001-6801-9603}$}
\affil{
$^{a}$College of Economic and Management Sciences, University of South Africa, Pretoria, 0003, South Africa \\
$^{b}$DSTI-NRF Centre of Excellence in Mathematical and Statistical Sciences (CoE-MaSS), South Africa
\newline \textit{E-mail}: \textup{jordakh@unisa.ac.za}\\
$^{c}$Department of Mathematics, Siksha `O' Anusandhan University, Bhubaneswar, 751030, Odisha, India
\newline \textit{E-mail}: \textup{pinakipk@gmail.com, pinakiprasadkar@soa.ac.in}
}

\begin{document}
\maketitle

\begin{abstract}
\noindent We establish quasi-orthogonality of arbitrary admissible order for $q$-Meixner polynomials \newline
$\{M_n(q^{-x};b,c;q)\}_{n=0}^{\infty}$ with respect to a positive discrete measure under simultaneous shifts of the parameters $b$ and $c$. This provides a construction not obtained by shifting either parameter alone. We then investigate the zero distribution and interlacing properties in the order-one case. In the parameter range considered, all zeros are real and simple, exactly one lies below $1$, and the remaining zeros lie above $q^{-1}$. Using connection formulas and $q$-contiguous relations, we establish interlacing results involving consecutive and nonconsecutive degrees, as well as related orthogonal $q$-Meixner polynomials obtained by one-step and two-step parameter shifts. In some cases where the original zero sets do not interlace, we identify explicit additional points that complete the interlacing.
\end{abstract}


\vspace{.2in} \noindent {\it {Keywords:}} Interlacing of zeros, Quasi-orthogonal polynomials, $q$-Meixner polynomials.

\section{Introduction}\label{Sec1}
Discrete orthogonal polynomials in the $q$-Askey scheme are basic hypergeometric families whose orthogonality
is expressed as a weighted sum over a $q$-linear or $q$-quadratic lattice rather than as an integral.
The $q$-Meixner polynomials belong to this class and are $q$-analogues of classical Meixner polynomials (cf. \cite[(9.10.1)]{Koekoek-Book}),
introduced by Josef Meixner in 1934 \cite{Meixner1934}.
Throughout this paper, we assume $0<q<1$.
The sequence of $q$-Meixner polynomials $\left\{M_n\left(q^{-x};b,c;q\right)\right\}_{n=0}^{\infty}$ defined by (see \cite[(14.13.1)]{Koekoek-Book})
\begin{equation*}
 M_n\left(q^{-x};b,c;q\right) =\ _2\phi_1\left(\begin{matrix} q^{-n}, q^{-x} \\ bq \end{matrix};q,-\frac{q^{n+1}}{c}\right),
\end{equation*}
is orthogonal with respect to the discrete weight
\begin{equation}\label{Mweight}
  w^q(x;b,c) = \frac{(bq;q)_x c^x q^{\binom{x}{2}}}{(q,-bcq;q)_x}
\end{equation}
when $0\leq bq<1$, $c>0$ and $x\in \{0,1,2,\dots\}.$ References to orthogonality (or quasi-orthogonality) with respect to the weight
$w^q(x;b,c)$ mean with respect to the discrete measure assigning mass $w^q(x;b,c)$ to the lattice point $z=q^{-x}$. The same convention applies to shifted parameter pairs. Since $0<q<1$, the transformation $z=q^{-x}$ is strictly increasing and maps the discrete support $x=0,1,2,\ldots$ onto the $q$-lattice $1,q^{-1},q^{-2},\ldots$. Moreover, the ordering of the zeros is preserved. It is well-known that the zeros of consecutive polynomials in an orthogonal sequence interlace. Furthermore, there can be at most one zero between two consecutive $q$-lattice points (cf. \cite[Theorem II.4.1]{Chihara-Book}, \cite[Theorem 3.41.2]{Szego-Book}).

Shifting the parameters of an orthogonal family beyond the ranges that ensure full orthogonality can give rise to quasi-orthogonal polynomials.
Originating from the classical work of Riesz \cite{Riesz},
a polynomial $Q_n(t)$ of degree $n$ is quasi-orthogonal of exact order $r$, $1\leq r<n$,
with respect to a positive discrete weight $w(t)$ on a support set $\{t_k\}_{k=0}^{\infty}$ if
\begin{equation*}
  \sum_{k=0}^{\infty} t_k^j Q_n(t_k) w(t_k) = 0 \quad \text{for } j = 0, 1, \dots, n - r - 1,
\end{equation*}
with $\sum_{k=0}^{\infty} t_k^{n-r} Q_n(t_k) w(t_k) \neq 0$.
Orthogonal polynomials with respect to a positive measure have real, simple zeros in the convex hull of the support, and the zeros of consecutive degrees strictly interlace. Quasi-orthogonality allows zeros outside the convex hull of the support, and the usual interlacing properties of orthogonal polynomials need not persist. This raises questions about both the location of the exceptional zeros and the separation of the remaining zeros.

In the classical setting, Driver and Jooste \cite{Driver2017} investigated the zeros of
order one quasi-orthogonal Meixner polynomials $M_n(x;\beta,c)$ arising from a single-parameter shift $-\beta,c\in(0,1)$,
where precisely one zero leaves the interval of orthogonality $(0,\infty)$.
They established interlacing properties within the quasi-orthogonal sequence
as well as with the nearest orthogonal counterparts $M_l(x;\beta+k,c)$ for $k\in\{1,2\}$.
In the $q$-Askey scheme, an analogous single-parameter approach was carried out by Kar and Gochhayat \cite{Kar2023-NA}
for order one quasi-orthogonal little $q$-Jacobi polynomials $p_n(z;a,b\mid q)$ with $aq^2,bq\in(0,1)$ and $aq>1$.
Using mixed $q$-contiguous relations, they analysed interlacing with the orthogonal sequences $p_m(z;aq^k,b\mid q)$ for $k=1,2$ and derived bounds for the exceptional least zero.

In 2018, Tcheutia et al. in \cite[Sec. 5]{Tcheutia2018} examined the parameter-shift properties of $q$-Meixner polynomials
and explained that shifting only the parameter $b$ does not give a suitable linear combination for quasi-orthogonality
with respect to a positive measure, whereas shifts in $c$ preserve full orthogonality.

In this paper, we show that simultaneous shifts of $b$ and $c$ yield quasi-orthogonality of arbitrary admissible order with respect to a positive discrete measure. We then concentrate on the order-one case, locating the exceptional zero and establishing interlacing relations involving different degrees and related orthogonal families. In cases where the original zero sets do not interlace, we identify explicit additional points that complete the interlacing.

For the analysis of zeros, we use the monic $q$-Meixner polynomial, defined by (see \cite[(14.13.4)]{Koekoek-Book})
\begin{equation}\label{s1e1}
  \widetilde{M}_{n}\left(z;b,c;q\right) := \frac{(-1)^n (b q;q)_n c^n}{q^{n^2}}\ _2\phi_1\left(\begin{matrix} q^{-n}, z \\ bq \end{matrix};q,-\frac{q^{n+1}}{c}\right).
\end{equation}
Note that, although the basic hypergeometric representation in \eqref{s1e1} may have singularities when $bq=q^{-m}$ for some nonnegative integer $m$, the monic polynomial $\widetilde M_n(z;b,c;q)$ remains well defined. Since the series terminates,
\[
\frac{(bq;q)_n}{(bq;q)_j}=(bq^{j+1};q)_{n-j},\, 0\leq j \leq n
\] and the quotient is interpreted after cancellation at singular parameter values.
Monic $q$-Meixner polynomials satisfy the three-term recurrence relation (see \cite[(14.13.3)]{Koekoek-Book})
\begin{equation}\label{s1e2}
  \widetilde{M}_{n}\left(z;b,c;q\right) = \left(z-A_n^q(b,c)\right) \widetilde{M}_{n-1}\left(z;b,c;q\right)
   + B_n^{q}(b,c) \widetilde{M}_{n-2}\left(z;b,c;q\right),
\end{equation}
where $A_n^q(b,c)=\frac{c \left(1-b q^n\right)+q \left(1-q^{n-1}\right) \left(c+q^{n-1}\right)+q^{2 n-1}}{q^{2 n-1}}$ and
$B_n^{q}(b,c) = \frac{c(q^{n-1}-1)(1-bq^{n-1})(c+q^{n-1})}{q^{4n-5}}.$

The remainder of this paper is organized as follows.
In Section \ref{Sec2}, we establish the quasi-orthogonality of $q$-Meixner polynomials with respect to a positive measure and locate the smallest zero in the order-one case.
Section \ref{Sec3} is devoted to deriving various interlacing properties of the zeros of quasi-orthogonal $q$-Meixner polynomials of order $1$,
including interlacing results among polynomials of different degrees and with their orthogonal counterparts via parameter shifts.

\section{Quasi-orthogonality of $q$-Meixner polynomials}\label{Sec2}
We begin this section by establishing quasi-orthogonality of arbitrary admissible order by expressing the polynomials with simultaneously shifted parameters as finite linear combinations of orthogonal $q$-Meixner polynomials. The order is determined with respect to the fixed weight $w^q(x;b,c)$.
\begin{theorem}\label{t1}
\noindent Let $n\geq 2$, $k=1,2,\ldots,n-1$, $c>0$, and $q^{k-1}<b<q^{-1}.$
Then $\widetilde{M}_n\left(z;\frac{b}{q^k},cq^k;q\right)$
is quasi-orthogonal of exact order $k$ with respect to the weight $w^q(x;b,c)$ and has at least $n-k$ zeros in $(1,\infty)$.
\end{theorem}
\begin{proof}
  Replacing $b$ by $b/q$ in \cite[(10d)]{Tcheutia2018-1}, we obtain the connection formula
  \begin{equation}\label{s2t1e1}
    \widetilde{M}_n\left(z;b/q, c q;q\right) = \widetilde{M}_n\left(z;b, c;q\right) +\displaystyle{\frac{c(1-q^n)}{q^{2n-1}}} \widetilde{M}_{n-1}\left(z;b, c;q\right).
  \end{equation}
Iterating this connection formula $k$ times gives
 \begin{equation}\label{s2t1e2}
    \widetilde{M}_{n}\left(z;b/q^k, c q^k;q\right) = \sum_{i=0}^{k} A_{n-i} \widetilde{M}_{n-i}\left(z;b,c;q\right),
  \end{equation}
where the $A_{n-i}$ are independent of $z$, $A_n=1$, and
$A_{n-k}=\displaystyle{\prod_{r=1}^{k}\frac{cq^{r-1}(1-q^{n-k+r})}{q^{2(n-k+r)-1}}\neq 0.}$
Since $q^{k-1}<b<q^{-1}$, we have $0<bq<1$, and hence the polynomials
$\{\widetilde{M}_m(z;b,c;q)\}_{m=0}^{\infty}$
form an orthogonal sequence with respect to the positive discrete $q$-Meixner measure with weight $w^q(x;b,c)$ defined in \eqref{Mweight}.
Moreover, $\frac{b}{q^k}q=bq^{1-k}>1$, so that the shifted parameter pair
lies outside the standard orthogonality range.
Evaluating \eqref{s2t1e2} at $z=q^{-x}$, multiplying the connection relation \eqref{s2t1e2} by
$w^q(x;b,c)(q^{-x})^j$
and summing over $x=0,1,2,\ldots$, we obtain
\begin{align*}
\sum_{x=0}^{\infty} w^q(x;b,c)(q^{-x})^j \widetilde{M}_n \left(q^{-x};\frac{b}{q^k},cq^k;q\right)
& = \sum_{x=0}^{\infty} w^q(x;b,c)(q^{-x})^j \sum_{i=0}^{k} A_{n-i}\widetilde{M}_{n-i} (q^{-x};b,c;q)\\
&=\sum_{i=0}^{k} A_{n-i} \sum_{x=0}^{\infty} w^q(x;b,c)(q^{-x})^j \widetilde{M}_{n-i} (q^{-x};b,c;q).
\end{align*}
If $j=0,1,\ldots,n-k-1$, then $j<n-i$ for every $i=0,1,\ldots,k$. By the orthogonality of
$\widetilde{M}_{n-i}(z;b,c;q),$
each term on the right-hand side vanishes. Therefore,
\[
\sum_{x=0}^{\infty} w^q(x;b,c)(q^{-x})^j \widetilde{M}_n\left(q^{-x};\frac{b}{q^k},cq^k;q\right) =0, \qquad j=0,1,\ldots,n-k-1.
\]
For $j=n-k$, all terms on the right-hand side vanish except the term corresponding to $i=k$.
Since $\widetilde{M}_{n-k}(z;b,c;q)$ is monic, $z^{n-k}-\widetilde M_{n-k}(z;b,c;q)$ has degree at most $n-k-1$.
Hence, by orthogonality,
\[ \sum_{x=0}^{\infty}w^q(x;b,c)(q^{-x})^{n-k} \widetilde M_{n-k}(q^{-x};b,c;q) = \sum_{x=0}^{\infty}w^q(x;b,c) \left[\widetilde M_{n-k}(q^{-x};b,c;q)\right]^2>0. \]
Since $A_{n-k}>0$, we have
\begin{equation*}
\sum_{x=0}^{\infty} w^q(x;b,c)(q^{-x})^j \widetilde{M}_n \left(q^{-x};\frac{b}{q^k},cq^k;q\right) =
A_{n-k} \sum_{x=0}^{\infty} w^q(x;b,c) \left[\widetilde{M}_{n-k}(q^{-x};b,c;q)\right]^2 >0.
\end{equation*}
Hence $\widetilde{M}_n\left(z;\frac{b}{q^k},cq^k;q\right)$
is quasi-orthogonal of exact order $k$ with respect to $w^q(x;b;c)$.
The assertion about its zeros follows from
\cite[Theorem 2]{Joulak2005}.
\end{proof}

\begin{remark}
The order of quasi-orthogonality depends on the reference weight. Although the polynomial in Theorem \ref{t1} has exact order $k$ with respect to $w^q(x;b,c)$, it may have lower order with respect to a weight obtained by an intermediate parameter shift.
Indeed, for $1\le r<k$,
\[ \widetilde M_n\!\left(z;\frac{b}{q^k},cq^k;q\right) = \widetilde M_n\!\left(z;\frac{\,b/q^{k-r}\,}{q^r}, (cq^{k-r})q^r;q \right). \]
Hence, whenever $q^{k-1}<b<q^{\,k-r-1}$, Theorem \ref{t1}, applied with $b$ replaced by $b/q^{k-r}$ and $c$ by $cq^{k-r}$,
shows that the same polynomial is quasi-orthogonal of exact order $r$ with respect to the corresponding intermediate $q$-Meixner measure.
In particular, when $q^{k-1}<b<1$, it is also quasi-orthogonal of order $k-1$ with respect to the measure with parameters $b/q$ and $cq$.
Thus, within the range of Theorem \ref{t1}, $1<b<q^{-1}$ is the subinterval in which none of these intermediate weights yields lower order of quasi-orthogonality.
\end{remark}
\begin{remark}\label{s2r2}
We now relate this construction to the classical Meixner family.
Using the standard (nonmonic) normalization $M_n$, let $0<c<1,$ $0<\beta<1,$ and $k=1,2,\dots,n-1.$
Substitute $q^{\beta-1}$ and $\displaystyle{\frac{c}{1-c}}$ for the first and second parameters, respectively, in the shifted $q$-Meixner polynomials, to obtain
\begin{eqnarray*}
{M}_{n}\left(q^{-x};\frac{q^{\beta-1}}{q^k}, \frac{c}{1-c} q^k;q\right) & = & {M}_{n}\left(q^{-x}; q^{\beta-k-1}, \frac{c}{1-c} q^k;q\right) \\
&=&\ _2\phi_1 \left(\begin{matrix} q^{-n}, q^{-x} \\ q^{\beta-k} \end{matrix}; q, -\displaystyle{ \frac{q^{n+1}}{\left(\frac{cq^k}{1-c}\right)}}\right) \\
 &=&\ _2\phi_1 \left(\begin{matrix} q^{-n}, q^{-x} \\ q^{\beta-k} \end{matrix}; q, -\displaystyle{\left(\frac{1-c}{c}\right) q^{n-k+1}}\right).
\end{eqnarray*}
Taking the limit $q\to1^-$, we obtain
\begin{equation*}
\lim\limits_{q\to1^-} {M}_{n}\left(q^{-x};\frac{q^{\beta-1}}{q^k}, \frac{c}{1-c} q^k;q\right) = M_n (x;\beta-k,c).
\end{equation*} Since $0<\beta < 1,$ the parameter $\beta - k$ of the limiting polynomial satisfies $-k < \beta - k < -k + 1,$
so $M_n\left(x;\beta - k, c\right)$ is quasi-orthogonal of order $k$.
This recovers the quasi-orthogonality of order $k$ described in \cite[Remark 4]{Jooste2013}.
\end{remark}

We now specialise to order one. Theorem \ref{t1} guarantees that at least $n-1$ zeros of $\widetilde{M}_n(z;b/q,cq;q)$ lie in $(1,\infty)$.
For $1<b<q^{-1}$, the following result shows that the remaining zero lies in $(-\infty,1)$.
\begin{theorem}\label{s2t2}
  Let $n\in\mathbb{N},$ $1<b<q^{-1}$ and $c>0.$ The smallest zero $z_{1,n}$  of $\widetilde{M}_n\left(z;\displaystyle{b/q}, cq; q\right)$, which is quasi-orthogonal of order $1$  with respect to the weight $w^q(x;b,c)$, satisfies \[z_{1,n}<1.\]
\end{theorem}
\begin{proof}
  Let $f_n(z)=\displaystyle{\frac{\widetilde{M}_{n}\left(z;b,c; q\right)}{\widetilde{M}_{n-1}\left(z;b,c; q\right)}}.$
  Using \eqref{s2t1e1} and $1<b<q^{-1}$, we obtain
  $$\displaystyle \frac{c(q^n-1)}{q^{2n-1}} < f_n(1)  = -\frac{c (1-bq^n)}{q^{2n-1}} < 0.$$
  Hence,  by \cite[Theorem 4(i)]{Joulak2005}, $z_{1,n} < 1.$
\end{proof}

\begin{remark}
The $q$-Meixner difference equation remains valid in the parameter range where the polynomials are quasi-orthogonal. In particular, let $P_n(z)=\widetilde{M}_n(z;b,c;q)$, where $0<q<1$, $c>0$, and $q^{-1}<b<q^{-2}$. Then, for $z\notin\{q^{-1},-bc/q\}$,
\begin{equation*}
\mathcal{D}_q^2 P_n(z) + \frac{c(1-bq)+q[1-z(1+q-q^n)]}{(1-q) (zq-1)(zq+bc)} \mathcal{D}_q P_n(z) + \frac{q(1-q^n)}{(1-q)^2 (zq-1) (zq+bc)} P_n(z) = 0,
\end{equation*}
where
\begin{equation*}
\mathcal{D}_q f(z) = \frac{f(qz)-f(z)}{z (q-1)},
\end{equation*}
with the value at $z=0$ defined by continuity for polynomials.
This is the standard $q$-Meixner $q$-difference equation \cite[(14.13.5)]{Koekoek-Book}, expressed using the $q$-derivative operator $D_q$.
After monic normalisation and clearing denominators, the left-hand side is a polynomial in $z$, $b$, and $c$,
and therefore extends directly from the orthogonality regime to the quasi-orthogonal parameter range.
\end{remark}

\section{Interlacing of zeros}\label{Sec3}
Throughout this section, we assume $q^{-1}<b<q^{-2}$ and $c>0$. Taking $k=1$ in Theorem \ref{t1} and replacing $b$ by $bq$ and $c$ by $c/q$, we see that, for $q^{-1}<b<q^{-2}$ and $c>0$, the polynomial $\widetilde{M}_n\left(z;b, c; q\right)$ is quasi-orthogonal of order $1$ with respect to the $q$-Meixner weight $w^q(x;bq,c/q)$. Moreover, under the same conditions, $\{\widetilde{M}_{l}\left(z;bq, c/q; q\right)\}_{l=0}^{\infty}$ and $\{\widetilde{M}_{l}\left(z;bq^2, c/q^2; q\right)\}_{l=0}^{\infty}$ are orthogonal $q$-Meixner polynomial sequences with respect to the weights $w^q(x;b q,c/q)$ and $w^q(x;b q^2,c/q^2)$ respectively.

In this section, we investigate interlacing of zeros within the quasi-orthogonal family and between this family and related families obtained by simultaneous parameter shifts.
We will denote the zeros of $\widetilde M_m(z;b,c;q)$, $\widetilde M_m(z;bq,c/q;q)$, and $\widetilde M_m(z;bq^2,c/q^2;q)$ by $z_{1,m}<\cdots<z_{m,m}$, $x_{1,m}<\cdots<x_{m,m}$, and $y_{1,m}<\cdots<y_{m,m}$, respectively throughout this section.

\subsection{Basic interlacing and zero location}
\begin{theorem}\label{s3t1}
  Let $n\ge 2,$ $n\in\mathbb{N},$ $q^{-1}<b<q^{-2},$ $c>0$. Then the zeros of $\widetilde{M}_{n}\left(z;b, c; q\right)$ interlace with
  those of $\widetilde{M}_{n}\left(z;bq, c/q; q\right)$ and $\widetilde{M}_{n-1}\left(z;bq, c/q; q\right)$ as follows:
  \begin{equation*}
    z_{1,n}<1<x_{1,n}<x_{1,n-1}<z_{2,n}<x_{2,n}<x_{2,n-1}<\cdots<z_{n-1,n}<x_{n-1,n}<x_{n-1,n-1}<z_{n,n}<x_{n,n},
  \end{equation*}
 where $\widetilde{M}_{n}\left(z;b, c; q\right)$ and $\widetilde{M}_{n}\left(z;bq, c/q; q\right)$,
  respectively, are the quasi-orthogonal (order $1$) and orthogonal $q$-Meixner polynomials,
  both with respect to the $q$-Meixner weight $w^q(x;bq,c/q)$.
\end{theorem}

\begin{proof}
  Replacing $b$ by $bq$ and $c$ by $c/q$ in \eqref{s2t1e1}, we obtain
  \begin{equation}\label{s3t1e1}
    \widetilde{M}_{n}\left(z;b, c; q\right) = \widetilde{M}_{n}\left(z;bq, c/q; q\right) + \displaystyle\frac{c (1-q^n)}{q^{2n}} \widetilde{M}_{n-1}\left(z;bq, c/q; q\right).
  \end{equation}
  Evaluating \eqref{s3t1e1} at the consecutive zeros $x_{i,n}$ and $x_{i+1,n}, i\in\{1,\ldots,n-1\},$ of $\widetilde{M}_{n}\left(z;bq, c/q; q\right),$ we get
  \begin{equation}\label{s3t1e2}
    \widetilde{M}_{n}\left(x_{i,n};b, c; q\right) \widetilde{M}_{n}\left(x_{i+1,n};b, c; q\right)
    = \displaystyle \frac{c^2 (1-q^n)^2}{q^{4n}} \widetilde{M}_{n-1}\left(x_{i,n};bq, c/q; q\right) \widetilde{M}_{n-1}\left(x_{i+1,n};bq, c/q; q\right).
  \end{equation}
  Since the zeros of the orthogonal polynomials $\widetilde{M}_{n}\left(z;bq, c/q; q\right)$ and $\widetilde{M}_{n-1}\left(z;bq, c/q; q\right)$ are interlacing,
  we conclude that the right-hand side of \eqref{s3t1e2} is negative and hence the zeros of
  $\widetilde{M}_{n}\left(z;b, c; q\right)$ and $\widetilde{M}_{n}\left(z;bq, c/q; q\right)$ are interlacing.
  Similarly, evaluating \eqref{s3t1e1} at the consecutive zeros $x_{i,n-1}$ and $x_{i+1,n-1}, i\in\{1,\ldots,n-2\},$ of $\widetilde{M}_{n-1}\left(z;bq, c/q; q\right),$ yields that
  the zeros of $\widetilde{M}_{n}\left(z;b, c; q\right)$ and \newline $\widetilde{M}_{n-1}\left(z;bq, c/q; q\right)$ interlace.
  In fact, from Theorem \ref{s2t2}, we have $z_{1,n}<1,$ and it is known that the orthogonal polynomials
  $\widetilde{M}_{n}\left(z;bq, c/q; q\right)$ and $\widetilde{M}_{n-1}\left(z;bq, c/q; q\right)$ have all their zeros in the interval $(1,\infty).$
\end{proof}

\begin{corollary}\label{s3c1}
  Under the hypotheses of Theorem \ref{s3t1}, the second smallest zero of quasi-orthogonal order 1 polynomial $\widetilde{M}_{n}\left(z;b, c; q\right)$
  satisfies the inequality $bq < q^{-1} < z_{2,n}.$
\end{corollary}
\begin{proof}
  Since $q^{-1} < b < q^{-2},$ we have $1 < bq < q^{-1}.$
  The orthogonal polynomial $\widetilde{M}_{n}\left(z;bq, c/q; q\right)$ has at most one zero in the interval $\left(1,q^{-1}\right),$
  because its orthogonality measure is supported on the $q$-lattice $\{1,q^{-1},q^{-2},\ldots\}.$
  Since a lattice point lies strictly between any two consecutive zeros of the orthogonal polynomial, $x_{2,n}>q^{-1}$. Theorem \ref{s3t1} therefore gives $z_{j,n}>q^{-1}$ for $j=3,\ldots,n$. Since $z_{1,n}<1$, only $z_{2,n}$ can lie in $(1,q^{-1})$.
  Evaluating $\widetilde{M}_{n}\left(z;b, c; q\right)$ at $z=1$ and $z=q^{-1}$ gives nonzero values of the same sign. Since all zeros are simple, the number of zeros in $(1,q^{-1})$ must therefore be even. Since $z_{2,n}$ is the only possible zero in this interval, there are none. Furthermore, $\widetilde M_n(q^{-1};b,c;q)\ne0$, so $z_{2,n}>q^{-1}$. Together with $bq<q^{-1}$, this proves the assertion.
\end{proof}

\subsection{Consecutive and non-consecutive degrees}
In this section, we examine the interlacing behavior of zeros of quasi-orthogonal $q$-Meixner polynomials of order $1$,
including both consecutive and non-consecutive degrees alongside their orthogonal counterparts.
Interlacing results involving consecutive degree polynomials cover all real zeros,
whereas those involving non-consecutive degrees are strictly confined to the interval $(1,\infty)$.

We begin by comparing consecutive members of the quasi-orthogonal family.
The following result shows that adjoining the endpoint $1$ to the zeros of the lower-degree polynomial completes the interlacing.

\begin{theorem}\label{s3t2}
Suppose that $n\ge 3,$ $n\in\mathbb{N},$ $q^{-1}<b<q^{-2},$ and $c>0.$
Then the zeros of $(z-1) \widetilde{M}_{n-1}\left(z;b, c; q\right)$ interlace with the zeros of $\widetilde{M}_{n}\left(z;b, c; q\right),$
where both $\widetilde M_{n-1}(z;b,c;q)$ and $\widetilde M_n(z;b,c;q)$ are quasi-orthogonal of order $1$ with respect to the $q$-Meixner weight $w^q(x;bq,c/q)$.
\end{theorem}

\begin{proof}
Using \cite[(14)]{Kar2020} and \cite[(8)]{Kar2025}, we obtain the following $q$-contiguous relation satisfied by the sequence of $q$-Meixner polynomials:
\begin{equation}\label{s3t2e1}
  \widetilde{M}_{n}\left(z;b, c; q\right) = (z-bq) \widetilde{M}_{n-1}\left(z;bq, c/q; q\right) - \displaystyle \frac{(c+q^n)(1-bq^n)}{q^{2n-1}} \widetilde{M}_{n-1}\left(z;b,c;q\right).
\end{equation}
Evaluating \eqref{s3t2e1} at the consecutive zeros $z_{i,n-1}$ and $z_{i+1,n-1},$ $i=1,\ldots,n-2,$
of $\widetilde{M}_{n-1}\left(z;b,c;q\right),$ gives
\begin{multline}\label{s3t2e2}
  \widetilde{M}_{n}\left(z_{i,n-1};b,c;q\right) \widetilde{M}_{n}\left(z_{i+1,n-1};b,c;q\right) \\
  = (z_{i,n-1}-bq) (z_{i+1,n-1}-bq) \widetilde{M}_{n-1}\left(z_{i,n-1};bq,c/q;q\right) \widetilde{M}_{n-1}\left(z_{i+1,n-1};bq,c/q;q\right).
\end{multline}
We know from Theorem \ref{s3t1} that the zeros of $\widetilde{M}_{n-1}\left(z;bq, c/q; q\right)$ and $\widetilde{M}_{n-1}\left(z;b, c; q\right)$ interlace.
By Corollary \ref{s3c1}, we have $1<bq<z_{2,n-1}.$
Then since $z_{1,n-1}<1$ implies $z_{1,n-1}-bq<0,$ the product $(z_{1,n-1}-bq) (z_{2,n-1}-bq)<0.$
However, $(z_{i,n-1}-bq) (z_{i+1,n-1}-bq) >0$ for $i\in\{2,\ldots,n-2\}.$
Hence from \eqref{s3t2e2}, $\widetilde{M}_{n}\left(z_{i,n-1};b, c; q\right)$ and $\widetilde{M}_{n}\left(z_{i+1,n-1};b, c; q\right)$
have the same sign for $i=1,$ and for $i=2,\ldots, n-2,$ they have different signs.
We also have from Theorem \ref{s3t1} that $z_{n,n}>x_{n-1,n-1}>z_{n-1,n-1}$ and $z_{1,n}<1.$
Therefore, we have exactly $(n-3)$ zeros of $\widetilde{M}_{n}\left(z;b, c; q\right)$ in between $z_{2,n-1}$ and $z_{n-1,n-1}.$
We still need to locate two more zeros that must be between $z_{1,n-1}$ and $z_{2,n-1}.$
Since $z_{1,n-1} , z_{1,n} < 1,$ the only possible configuration of zeros is
\begin{equation}\label{s3t2e3}
  z_{1,n-1} < z_{1,n} < 1 < z_{2,n} < z_{2,n-1} < \cdots < z_{n-1,n-1} < z_{n,n},
\end{equation}
which shows the required interlacing result.
\end{proof}

We now consider polynomials whose degrees differ by two and show that one of the recurrence coefficients supplies an additional interlacing point.
\begin{theorem}\label{s3t3}
Assume that $n\ge 3,$ $n\in\mathbb{N},$ $q^{-1}<b<q^{-2},$ and $c>0$
and that the quasi-orthogonal order $1$ polynomials $\widetilde{M}_{n}\left(z;b, c; q\right)$
and $\widetilde{M}_{n-2}\left(z;b, c; q\right)$ with respect to the $q$-Meixner weight $w^q(x;bq,c/q)$ have no common zeros.
Then the zeros of $\left(z-A_n^q(b,c)\right)$ $\widetilde{M}_{n-2}\left(z;b, c; q\right)$
interlace with those of $\widetilde{M}_{n}\left(z;b, c; q\right)$ in the interval $(1,\infty).$
\end{theorem}

\begin{proof}
  From \eqref{s1e2}, it is clear that $z_{i,n}\neq A_n^q(b,c),$ $i\in\{1,\ldots,n\},$
  since the polynomials $\widetilde{M}_{n}\left(z;b, c; q\right)$ and $\widetilde{M}_{n-2}\left(z;b, c; q\right)$ have no zeros in common.
  Evaluating \eqref{s1e2} at the consecutive zeros $z_{i,n}$ and $z_{i+1,n},$ $i\in\{2,\ldots,n-1\},$ of $\widetilde{M}_{n}\left(z;b, c; q\right)$ yields
  \begin{multline}\label{s3t3e1}
  \left(z_{i,n}-A_n^q(b,c)\right) \left(z_{i+1,n}-A_n^q(b,c)\right) \widetilde{M}_{n-1}\left(z_{i,n};b, c; q\right) \widetilde{M}_{n-1}\left(z_{i+1,n};b, c; q\right)\\
  = \left(B_n^q (b,c)\right)^2 \widetilde{M}_{n-2}\left(z_{i,n};b, c; q\right) \widetilde{M}_{n-2}\left(z_{i+1,n};b, c; q\right).
  \end{multline}
  We know from Theorem \ref{s3t2} that the zeros, that are greater than $1,$
  of $\widetilde{M}_{n-1}\left(z;b, c; q\right)$ interlace with those of  $\widetilde{M}_{n}\left(z;b, c; q\right).$
  Also, the LHS of \eqref{s3t3e1} is negative,
  if $A_n^q(b,c)\notin\left(z_{i,n},z_{i+1,n}\right),$ $i\in\{2,\ldots,n-1\}.$
  That means the LHS of \eqref{s3t3e1} is positive when $A_n^q(b,c)\in\left(z_{j,n},z_{j+1,n}\right)$ for one $j\in\{2,\ldots,n-1\}.$
  That is, $\widetilde{M}_{n-2}\left(z_{i,n};b, c; q\right) \widetilde{M}_{n-2}\left(z_{i+1,n};b, c; q\right)<0$ for each $i\in\{2,\ldots,n-1\},$
  except when $A_n^q(b,c)\in\left(z_{j,n},z_{j+1,n}\right)$ for one $j\in\{2,\ldots,n-1\}.$
  Since $\widetilde{M}_{n}\left(z;b, c; q\right)$ has $n-1$ zeros in $(1,\infty),$ it has $n-2$ intervals with endpoints at its consecutive zeros that are $>1.$
  Also, $\widetilde{M}_{n-2}\left(z;b, c; q\right)$ has exactly $n-3$ zeros in $(1,\infty).$
  Therefore, we deduce that the zeros of $\widetilde{M}_{n-2}\left(z;b, c; q\right)$
  and the point $A_n^q(b,c)$ interlace with those of $\widetilde{M}_{n}\left(z;b, c; q\right)$ in the interval $(1,\infty).$
\end{proof}

A related comparison with the orthogonal polynomial of degree $n-2$ involves the explicitly determined point $E_n^q(b,c)$.
\begin{theorem}\label{s3t4}
For $n\ge 3,$ $n\in\mathbb{N},$ $q^{-1}<b<q^{-2},$ and $c>0,$ the zeros of
\[
(z-E_n^q(b,c))\widetilde M_{n-2}(z;bq,c/q;q),
\qquad
E_n^q(b,c)=\frac{c+q^n-bcq^n}{q^{2n-1}},
\]
where $\widetilde M_{n-2}(z;bq,c/q;q)$ is orthogonal, interlace in $(1,\infty)$ with those of the quasi-orthogonal order $1$ polynomial $\widetilde M_n(z;b,c;q),$
both with respect to the $q$-Meixner weight $w^q(x;bq,c/q)$, provided the two polynomials have no common zeros.
\end{theorem}

\begin{proof}
Using \eqref{s1e2} in \eqref{s3t2e1}, we obtain
\begin{equation}\label{s3t4e1}
\widetilde{M}_{n}\left(z;b, c; q\right) = \displaystyle \left(z-E_n^q(b,c)\right) \widetilde{M}_{n-1}\left(z;b, c; q\right) - \frac{c \left(1-q^{n-1}\right)}{q^{2 n-2}} (z-bq) \widetilde{M}_{n-2}\left(z;bq, c/q; q\right).
\end{equation}
It is clear from \eqref{s3t4e1} that $z_{i,n}\ne E_n^q (b,c),$ $i\in\{1,\ldots,n\},$
if not, this would contradict our assumption that the polynomials $\widetilde{M}_{n}\left(z;b, c; q\right)$ and $\widetilde{M}_{n-2}\left(z;bq, c/q; q\right)$ have no zeros in common.
Evaluating \eqref{s3t4e1} at the consecutive zeros $z_{i,n}$ and $z_{i+1,n},$ $i\in\{2,\ldots,n-1\},$ of $\widetilde{M}_{n}\left(z;b, c; q\right),$ we get
\begin{multline}\label{s3t4e2}
  \left(z_{i,n}-E_n^q(b,c)\right) \left(z_{i+1,n}-E_n^q(b,c)\right) \widetilde{M}_{n-1}\left(z_{i,n};b, c; q\right) \widetilde{M}_{n-1}\left(z_{i+1,n};b, c; q\right)\\
   = \displaystyle \frac{c^2 \left(1-q^{n-1}\right)^2}{q^{4n-4}} (z_{i,n}-bq) (z_{i+1,n}-bq) \widetilde{M}_{n-2}\left(z_{i,n};bq, c/q; q\right) \widetilde{M}_{n-2}\left(z_{i+1,n};bq, c/q; q\right).
\end{multline}
Since $bq < z_{2,n}$ by Corollary \ref{s3c1}, $(z_{i,n}-bq) (z_{i+1,n}-bq)>0$ for any $i\in\{2,\ldots,n-1\}.$
Also by Theorem \ref{s3t2}, the zeros of $\widetilde{M}_{n}\left(z;b, c; q\right)$ interlace with those of $\widetilde{M}_{n-1}\left(z;b, c; q\right)$ on $(1,\infty).$
Hence the left-hand side of \eqref{s3t4e2} is negative if $E_n^q(b,c)\notin (z_{i,n},z_{i+1,n})$ for $i\in\{2,\ldots,n-1\}.$
If $E_n^q(b,c)\notin (z_{i,n},z_{i+1,n})$ for any $i\in\{2,\ldots,n-1\},$
then from \eqref{s3t4e2}, we deduce that the zeros of
$\widetilde{M}_{n-2}\left(z;bq, c/q; q\right)$ and $\widetilde{M}_{n}\left(z;b, c; q\right)$ interlace in the interval $(1,\infty).$
In addition, since $E_n^q(b,c)\notin(z_{2,n},z_{n,n}),$
the zeros of $(z-E_n^q(b,c)) \widetilde{M}_{n-2}\left(z;bq, c/q; q\right)$
and the zeros, that are $>1,$ of $\widetilde{M}_{n}\left(z;b, c; q\right)$ interlace.
On the other hand, if $E_n^q(b,c)\in(z_{i,n},z_{i+1,n}),$ $i\in\{2,\ldots,n-1\},$
then $E_n^q(b,c)$ must lie in exactly one interval $(z_{j,n},z_{j+1,n}),$ where $j\in\{2,\ldots,n-1\}.$
Therefore, $\widetilde{M}_{n-2}\left(z_{j,n};bq, c/q; q\right)\widetilde{M}_{n-2}\left(z_{j+1,n};bq, c/q; q\right)>0$ for one $j\in\{2,\ldots,n-1\},$
however, $\widetilde{M}_{n-2}\left(z_{i,n};bq, c/q; q\right)\widetilde{M}_{n-2}\left(z_{i+1,n};bq, c/q; q\right)<0$ for all other $i\in\{2,\ldots,n-1\}.$
Finally, for $i=1,$ \eqref{s3t4e2} becomes
\begin{multline}\label{s3t4e3}
  \left(z_{1,n}-E_n^q(b,c)\right) \left(z_{2,n}-E_n^q(b,c)\right) \widetilde{M}_{n-1}\left(z_{1,n};b, c; q\right) \widetilde{M}_{n-1}\left(z_{2,n};b, c; q\right)\\
   = \displaystyle \frac{c^2 \left(1-q^{n-1}\right)^2}{q^{4n-4}} (z_{1,n}-bq) (z_{2,n}-bq) \widetilde{M}_{n-2}\left(z_{1,n};bq, c/q; q\right) \widetilde{M}_{n-2}\left(z_{2,n};bq, c/q; q\right).
\end{multline}
Since $E_n^q(b,c)\in(z_{j,n},z_{j+1,n})$ where $2\le j \le n-1,$ $E_n^q\notin (z_{1,n},z_{2,n}).$
Also, $z_{1,n} < 1 < b q < z_{2,n}$ and $\widetilde{M}_{n-1}\left(z_{1,n};b, c; q\right) \widetilde{M}_{n-1}\left(z_{2,n};b, c; q\right) > 0$ by \eqref{s3t2e3}.
Thus, we deduce from \eqref{s3t4e3} that there is exactly one zero of $\widetilde{M}_{n-2}\left(z;bq, c/q; q\right)$ in $(z_{1,n},z_{2,n}).$
Therefore, in both of the cases, we have the zeros of $\widetilde{M}_{n-2}\left(z;bq, c/q; q\right)$
and the point $E_n^q(b,c)$ interlace with the $(n-1)$ zeros of $\widetilde{M}_{n}\left(z;b, c; q\right),$
that are $>1,$ provided $\widetilde{M}_{n}\left(z;b, c; q\right)$ and $\widetilde{M}_{n-2}\left(z;bq, c/q; q\right)$ do not have common zeros. 
\end{proof}

\subsection{An additional parameter shift}
We now turn to the parameter pair $(bq^2,c/q^2)$.
We first establish the relative ordering of the zeros of the two shifted orthogonal families
$\widetilde{M}_{n}\left(z;bq, c/q; q\right)$ and $\widetilde{M}_{n}\left(z;bq^2, c/q^2; q\right)$.

Combining \eqref{s3t1e1} and \eqref{s3t2e1} yields
  \begin{multline}\label{s3l1e1}
    \frac{(1-bq^{n+1})(c+q^{n+1})}{q^{2n+1}} \widetilde{M}_{n}\left(z;b, c; q\right) = \left(z-F_n^q(b,c)\right) \widetilde{M}_{n}\left(z;bq, c/q; q\right)\\
   - \left(z-bq^2\right) \widetilde{M}_{n}\left(z;bq^2, c/q^2; q\right),
  \end{multline}
  where $F_n^q(b,c) = \displaystyle \frac{b q^{n+2} (q+1)-(q + c (1-b q))}{q^{n+1}}.$

\begin{lemma}\label{s3l1}
  Suppose that $n\in\mathbb{N},$ $q^{-1}<b<q^{-2},$ and $c>0.$
  Then the zeros of the $q$-Meixner polynomial $\widetilde{M}_{n}\left(z;bq, c/q; q\right),$
  orthogonal with respect to the $q$-Meixner weight $w^q(x;bq,c/q)$,
  interlace with the zeros
  of the $q$-Meixner polynomial $\widetilde{M}_{n}\left(z;bq^2, c/q^2; q\right),$
  orthogonal with respect to the $q$-Meixner weight $w^q(x;bq^2,c/q^2)$.
  The zeros interlace in the order
\begin{equation}\label{s3l1e3}
x_{1,n}<y_{1,n}<x_{2,n}<y_{2,n}<\cdots<x_{n,n}<y_{n,n}.
\end{equation}
\end{lemma}
\begin{proof}
  First suppose that $n\geq2$.
In addition to the two orthogonal polynomials in the statement,
we use $\widetilde M_n(z;b,c;q)$, which is quasi-orthogonal of
order~$1$ with respect to the weight $w^q(x;bq,c/q)$.
Thus Theorem~\ref{s3t1} applies to $\widetilde M_n(z;b,c;q)$ and
$\widetilde M_n(z;bq,c/q;q)$ with respect to this same measure.

  Evaluating \eqref{s3l1e1} at consecutive zeros $x_{i,n}$ and $x_{i+1,n},$ $i\in\{1,\ldots,n-1\},$
  of $\widetilde{M}_{n}\left(z;bq, c/q; q\right),$ yields
  \begin{multline}\label{s3l1e2}
    \frac{(1-bq^{n+1})^2 (c+q^{n+1})^2}{q^{4n+2}} \widetilde{M}_{n}\left(x_{i,n};b, c; q\right) \widetilde{M}_{n}\left(x_{i+1,n};b, c; q\right) \\
    = \left(x_{i,n}-bq^2\right) \left(x_{i+1,n}-bq^2\right) \widetilde{M}_{n}\left(x_{i,n};bq^2, c/q^2; q\right) \widetilde{M}_{n}\left(x_{i+1,n};bq^2, c/q^2; q\right).
  \end{multline}
  From Theorem \ref{s3t1}, the zeros of $\widetilde{M}_{n}\left(z;b, c; q\right)$ interlace with those of $\widetilde{M}_{n}\left(z;bq, c/q; q\right),$
  which implies the left-hand side of \eqref{s3l1e2} is strictly negative.
  Since $q^{-1}<b<q^{-2},$ $0<bq^2<1$ and $x_{i,n}>1$ for every $i,$ the product of first two factors on the right-hand side of \eqref{s3l1e2} is positive.
  Hence $\widetilde M_n(z;bq^2,c/q^2;q)$ has a zero in each interval
$(x_{i,n},x_{i+1,n})$, $i=1,\ldots,n-1$. It remains to locate the final zero. Let $z_{n,n}$ be the largest zero
of $\widetilde M_n(z;b,c;q)$. Theorem \ref{s3t1}  gives $z_{n,n}<x_{n,n}$.
Since $\widetilde M_n(z;b,c;q)$ is monic and all its zeros are real,
\[
\widetilde M_n(x_{n,n};b,c;q)>0.
\] Evaluating \eqref{s3l1e1}  at $x_{n,n}$ therefore gives
\[
\widetilde M_n(x_{n,n};bq^2,c/q^2;q)
=-\frac{(1-bq^{n+1})(c+q^{n+1})}
{q^{2n+1}(x_{n,n}-bq^2)}
\widetilde M_n(x_{n,n};b,c;q)<0,
\]
where $1-bq^{n+1}>0$ and $x_{n,n}-bq^2>0$.
As $\widetilde M_n(z;bq^2,c/q^2;q)$ is monic, it is positive for
sufficiently large $z$, and so it has a zero in $(x_{n,n},\infty)$.
The $n$ disjoint intervals just identified contain all its $n$ zeros.
It follows that the zeros of $\widetilde{M}_{n}\left(z;bq, c/q; q\right)$ and $\widetilde{M}_{n}\left(z;bq^2, c/q^2; q\right)$ interlace as
  \begin{equation*}
  x_{1,n} < y_{1,n} < x_{2,n} < y_{2,n} < \cdots < x_{n,n} < y_{n,n}.
  \end{equation*}
  For $n=1$, the formula
\[
\widetilde M_1(z;b,c;q)=z-1-\frac{c(1-bq)}{q}
\]
gives
\[
x_{1,1}=1+\frac{c(1-bq^2)}{q^2},\qquad
y_{1,1}=1+\frac{c(1-bq^3)}{q^3}.
\]
Thus $y_{1,1}-x_{1,1}=c(1-q)/q^3>0$, completing the proof.
\end{proof}
\begin{theorem}\label{s3t5}
  Let $n\ge 2,$ $n\in\mathbb{N},$ $q^{-1}<b<q^{-2},$ and $c>0.$
  Assume that the polynomial $\widetilde{M}_{n}\left(z; b, c; q\right)$,
  which is quasi-orthogonal of order $1$ with respect to the $q$-Meixner weight $w^q(x; bq, c/q)$,
  and the polynomial $\widetilde{M}_{n}\left(z; bq^2, c/q^2; q\right)$,
  orthogonal with respect to the $q$-Meixner weight
  $w^q(x; bq^2, c/q^2)$, are co-prime.
  Then the zeros of $\widetilde{M}_{n}\left(z;b, c; q\right)$ and $\widetilde{M}_{n}\left(z;bq^2, c/q^2; q\right)$ interlace,
  if $F_n^q(b,c) \notin \left(y_{i,n},y_{i+1,n}\right)$ for every $i\in\{1,\ldots,n-1\},$ where $F_n^q(b,c)$ is as given in \eqref{s3l1e1}.
  Otherwise, the zeros of \newline $\left(z-F_n^q(b,c)\right) \widetilde{M}_{n}\left(z;b, c; q\right)$ and $(z-1) \widetilde{M}_{n}\left(z;bq^2, c/q^2; q\right)$ interlace.
\end{theorem}
\begin{proof}
  Since $\widetilde{M}_{n}\left(z;b, c; q\right)$ and $\widetilde{M}_{n}\left(z;bq^2, c/q^2; q\right)$
  have no common zeros, \eqref{s3l1e1} implies that $F_n^q(b,c) \ne y_{i,n}$ for every $i=1,\ldots,n.$
  Indeed, if $F_n^q(b,c) = y_{i,n}$, then both terms on the right-hand side of \eqref{s3l1e1}
  vanish at $z=F_n^q(b,c)$ giving $\widetilde{M}_{n}\left(F_n^q(b,c);b, c; q\right)=0$.
  Thus $F_n^q(b,c)$ would be a common zero of the two polynomials, a contradiction.
  Evaluating \eqref{s3l1e1} at the consecutive zeros
  $y_{i,n}$ and $y_{i+1,n},$ $i\in\{1,\ldots,n-1\},$ of $\widetilde{M}_{n}\left(z;bq^2, c/q^2; q\right),$ we obtain
  \begin{multline}\label{s3t5e1}
  \frac{(1-bq^{n+1})^2 (c+q^{n+1})^2}{q^{4n+2}} \widetilde{M}_{n}\left(y_{i,n};b, c; q\right) \widetilde{M}_{n}\left(y_{i+1,n};b, c; q\right) \\
    = \left(y_{i,n}-F_n^q(b,c)\right) \left(y_{i+1,n}-F_n^q(b,c)\right) \widetilde{M}_{n}\left(y_{i,n};bq, c/q; q\right) \widetilde{M}_{n}\left(y_{i+1,n};bq, c/q; q\right).
  \end{multline}

  Since $F_n^q(b,c) \notin \left(y_{i,n}, y_{i+1,n}\right)$ for each $i \in \{1, \dots, n-1\},$
  the factors $(y_{i,n} - F_n^q(b,c))$ and $(y_{i+1,n} - F_n^q(b,c))$ share the same sign, making their product strictly positive.
  Furthermore, Lemma \ref{s3l1} implies that $\widetilde{M}_{n}\left(y_{i,n};bq, c/q; q\right) \widetilde{M}_{n}\left(y_{i+1,n};bq, c/q; q\right) < 0.$
  Because every zero of $\widetilde{M}_{n}\left(z; bq^2, c/q^2; q\right)$ is $> 1,$
  \eqref{s3t5e1} establishes that these $n$ zeros interlace with the $n-1$ zeros of $\widetilde{M}_{n}\left(z; b, c; q\right)$ on $(1, \infty).$
  Combined with the fact that the smallest zero of $\widetilde{M}_{n}\left(z; b, c; q\right)$ is strictly less than $1,$
  we conclude that the zeros must be ordered as $z_{i,n} < y_{i,n} < z_{i+1,n}$ for $i \in \{1, \dots, n-1\}$ and $z_{n,n} < y_{n,n}.$

  Now consider the complementary case, where $F_n^q(b,c) \in \left(y_{i,n},y_{i+1,n}\right),$ $i\in\{1,\ldots,n-1\}.$
  Then it necessarily belongs to exactly one interval $\left(y_{j,n},y_{j+1,n}\right),$ for one $j\in\{1, \ldots, n-1\}$
  and hence from \eqref{s3t5e1}, $\widetilde{M}_{n}\left(y_{j,n};b, c; q\right) \widetilde{M}_{n}\left(y_{j+1,n};b, c; q\right) > 0.$
  This implies the exceptional interval $(y_{j,n},y_{j+1,n})$ contains an even number of zeros of $\widetilde{M}_{n}\left(z;b, c; q\right)$.

  For $i\neq j$, the product $\widetilde{M}_{n}\left(y_{i,n};b, c; q\right) \widetilde{M}_{n}\left(y_{i+1,n};b, c; q\right)<0,$
  so each of the $n-2$ non-exceptional intervals contains an odd number of zeros of $\widetilde{M}_{n}\left(z;b, c; q\right)$,
  and therefore at least one.
  Since the polynomial $\widetilde{M}_{n}\left(z;b, c; q\right)$ has $n-1$ simple zeros in $(1,\infty)$,
  only one zero remains after placing one in each of the $n-2$ non-exceptional intervals.
  Consequently, none of these intervals can contain additional zeros, and the exceptional interval contains none.
  If the exceptional interval contains two or more zeros, then the total number of zeros would exceed $n-1$.
  The remaining one zero must lie in $(1,y_{1,n})$ or $(y_{n,n},\infty)$.
  By Theorem \ref{s3t1} and Lemma \ref{s3l1}, $z_{n,n}<x_{n,n}<y_{n,n}$,
  which excludes the latter possibility.
  Hence we have $z_{1,n}<1<z_{2,n}<y_{1,n}$.
  Now $(z-F_n^q(b,c))\widetilde{M}_{n}\left(z; b, c; q\right)$ has one zero $<1$,
  one zero in $(1,y_{1,n})$, exactly one zero in each interval $(y_{i,n},y_{i+1,n})$ for every $i\in\{1,\ldots,n-1\}$,
  and the exceptional interval contains $F_n^q (b,c)$.
  These account for all its $n+1$ zeros.
  They therefore strictly interlace with the zeros $1,\ y_{1,n},\ldots,y_{n,n}$,
  of $(z-1)\widetilde{M}_{n}\left(z; bq^2, c/q^2; q\right)$.
\end{proof}

In Theorem \ref{s3t5}, we compared the zeros of the degree $n$ quasi-orthogonal polynomial $\widetilde{M}_{n}\left(z; b, c; q\right)$ with the zeros of the degree $n$ orthogonal polynomial $\widetilde{M}_{n}\left(z; bq^2, c/q^2; q\right)$. Our next result compares the zeros of the degree $n$ quasi-orthogonal polynomial $\widetilde{M}_{n}\left(z; b, c; q\right)$ with the zeros of the degree $n-1$ orthogonal polynomial $\widetilde{M}_{n-1}\left(z; bq^2, c/q^2; q\right)$.
\begin{theorem}\label{s3t6}
Let $n \ge 2$, $n \in \mathbb{N}$, $q^{-1} < b < q^{-2}$, $c > 0$
and assume $q(1-bq^{n+1})+c(1-bq)\ne 0,  \text{ equivalently, }  c\ne c_*:=\displaystyle\frac{q(1-bq^{n+1})}{bq-1}$.
Then the zeros of the polynomial $(z - bq^2) \widetilde{M}_{n-1}\left(z; bq^2, c/q^2; q\right)$,
where $\widetilde{M}_{n-1}\left(z; bq^2, c/q^2; q\right)$ is orthogonal with respect to
the $q$-Meixner weight $w^q(x; bq^2, c/q^2)$,
and the polynomial $\widetilde{M}_{n}\left(z; b, c; q\right)$,
which is quasi-orthogonal of order $1$ with respect to the $q$-Meixner weight
$w^q(x; bq, c/q)$, interlace.
In addition, if $F_n^q(b,c) < z_{1,n},$
then the zeros of $\widetilde{M}_{n}\left(z;b, c; q\right)$
and $\widetilde{M}_{n-1}\left(z;bq^2, c/q^2; q\right)$ interlace,
where $F_n^q(b,c)$ is defined immediately after \eqref{s3l1e1}.
\end{theorem}
\begin{proof}
   The connection relation involving the polynomials $\widetilde{M}_{n}\left(z;b, c; q\right)$,
   $\widetilde{M}_{n}\left(z;bq, c/q; q\right)$,\\ and $\widetilde{M}_{n-1}\left(z;bq^2, c/q^2; q\right)$ is given by
   \begin{multline}\label{s3t6e1-1}
   \frac{\left(1-b q^{n+1}\right) \left(c+q^{n+1}\right)}{c \left(1-q^n\right)} \widetilde{M}_{n}\left(z;b, c; q\right)
   = \frac{q^{n+1} \left(1-b q^{n+1}\right) + c q^n (1-b q)}{c \left(1-q^n\right)} \widetilde{M}_{n}\left(z;bq, c/q; q\right) \\
   + (z-bq^2) \widetilde{M}_{n-1}\left(z;bq^2, c/q^2; q\right).
  \end{multline}
  Here, $\widetilde{M}_{n}\left(z;b, c; q\right)$ is quasi-orthogonal of order~$1$
  and $\widetilde{M}_{n}\left(z;bq, c/q; q\right)$ is orthogonal with respect to the discrete measure having masses $w^q(x;bq,c/q)$,
  whereas $\widetilde{M}_{n-1}\left(z;bq^2, c/q^2; q\right)$ is orthogonal with respect to the discrete measure
  having masses $w^q(x;bq^2,c/q^2)$. Both measures are supported at $z=q^{-x}$, $x=0,1,2,\ldots$.
  Evaluating \eqref{s3t6e1-1} at the consecutive zeros $z_{i,n}$ and $z_{i+1,n},$ $i\in\{1, \ldots , n-1\},$ of $\widetilde{M}_{n}\left(z;b, c; q\right),$ gives
  \begin{multline}\label{s3t6e2-1}
    \frac{q^{2n} \left[q(1-bq^{n+1})+c(1-b q)\right]^2}{c^2 (1-q^n)^2} \widetilde{M}_{n}\left(z_{i,n};bq, c/q; q\right) \widetilde{M}_{n}\left(z_{i+1,n};bq, c/q; q\right) \\
    = (z_{i,n}-bq^2) (z_{i+1,n}-bq^2) \widetilde{M}_{n-1}\left(z_{i,n};bq^2, c/q^2; q\right) \widetilde{M}_{n-1}\left(z_{i+1,n};bq^2, c/q^2; q\right).
  \end{multline}
  By Theorem \ref{s3t1}, the left-hand side of \eqref{s3t6e2-1} is negative for every $i\in\{1,\ldots,n-1\}.$
  Since $q^{-1}<b<q^{-2}$ implies $0<bq^2<1$ and $z_{2,n},z_{3,n}, \ldots , z_{n,n} > 1,$
  the product $(z_{i,n}-bq^2) (z_{i+1,n}-bq^2) >0$ for every $i\in\{2, \ldots , n-1\}.$
  It then follows from \eqref{s3t6e2-1} that
  $\widetilde{M}_{n-1}\left(z_{i,n};bq^2, c/q^2; q\right) \widetilde{M}_{n-1}\left(z_{i+1,n};bq^2, c/q^2; q\right)<0$
  for every $i\in\{2, \ldots , n-1\}.$
  Consequently, we have located $n-2$  zeros of $\widetilde{M}_{n}\left(z;bq^2, \blue{c/q^2}; q\right)$, one in each interval $(z_{i,n},z_{i+1,n})$, $i=2,\dots,n-1$.

  To locate $bq^2$ relative to $z_{1,n}$, we evaluate \eqref{s3t6e1-1} at $z_{1,n}$ to get
  \begin{equation}\label{s3t6e3-1}
    (z_{1,n}-bq^2) \widetilde{M}_{n-1}\left(z_{1,n};bq^2, c/q^2; q\right) = - \frac{q^n (bq-1) (c_* - c)}{c (1-q^n)}\widetilde{M}_{n}\left(z_{1,n};bq, c/q; q\right).
  \end{equation}
  Since $z_{1,n}$ lies strictly to the left of all zeros of the monic polynomials
  $\widetilde{M}_{n}\left(z;bq, c/q; q\right)$ and\newline $\widetilde{M}_{n-1}\left(z;bq^2, c/q^2; q\right)$,
  their signs at $z = z_{1,n}$ are given by
\begin{equation*}
  \operatorname{sgn}\left(\widetilde{M}_{n}\left(z_{1,n};bq, c/q; q\right)\right) = (-1)^n
  \quad \text{and} \quad
  \operatorname{sgn}\left(\widetilde{M}_{n-1}\left(z_{1,n};bq^2, c/q^2; q\right)\right) = (-1)^{n-1}.
  \end{equation*}
  Comparing signs on both sides of \eqref{s3t6e3-1} leads to
  \begin{equation*}
  \operatorname{sgn}\left(z_{1,n} - bq^2\right) = \operatorname{sgn}(c_* - c).
  \end{equation*}
  As a consequence:
  \begin{itemize}
    \item If $c < c_*$, then $bq^2 < z_{1,n} < 1$, and \eqref{s3t6e2-1} implies that
            $\widetilde{M}_{n-1}\left(z; bq^2, c/q^2; q\right)$ has a zero in $(z_{1,n}, z_{2,n})$, and hence in $(1, z_{2,n})$.
    \item If $c > c_*$, then $z_{1,n} < bq^2 < 1$, and \eqref{s3t6e2-1} shows that
            $\widetilde{M}_{n-1}\left(z; bq^2, c/q^2; q\right)$ has the same non-zero sign at $z_{1,n}$ and $z_{2,n}$. Since its zeros are simple, their number in $(z_{1,n},z_{2,n})$ is even. The $n-2$ intervals $(z_{i,n},z_{i+1,n})$, $i=2,\ldots,n-1$, each contain an odd number of zeros and therefore at least one. As the polynomial has degree $n-1$, the first interval contains no zeros, and each of the other intervals contains exactly one. Its remaining zero lies in $(z_{n,n},\infty)$, because all its zeros exceed $1>z_{1,n}$.
    \end{itemize}
In either case, the zeros of $(z-bq^2) \widetilde{M}_{n-1}\left(z; bq^2, c/q^2; q\right)$ interlace with those of $\widetilde{M}_{n}\left(z; b, c; q\right)$,
completing the proof of the first assertion.

  For the second assertion, the definition of $F_{n}^{q}(b,c)$ gives
  \[F_n^q(b,c) = \frac{bq^{n+2}(q+1)-\bigl(q+c(1-bq)\bigr)}{q^{n+1}} = bq^2 - \frac{(bq - 1) (c_* - c)}{q^{n+1}}.\]
  The preceding sign analysis reveals that
    \[
       \begin{array}{ll}
        c < c_* :& F_n^q(b,c) < bq^2 < z_{1,n},\\[2pt]
        c_* < c :& z_{1,n}    < bq^2 < F_n^q(b,c).
       \end{array}
    \]
Therefore, $F_n^q(b,c)<z_{1,n}$ holds if and only if $c < c_*$.
Under this condition, all $n-1$ zeros of $\widetilde{M}_{n-1}\left(z;bq^2, c/q^2; q\right)$
lie between consecutive zeros of $\widetilde{M}_{n}\left(z;b, c; q\right)$.
Writing them as $y_{1,n-1}<\cdots<y_{n-1,n-1}$, we obtain
\[
z_{i,n}<y_{i,n-1}<z_{i+1,n},\qquad i=1,\ldots,n-1,
\]
which proves the second assertion.
\end{proof}
\begin{remark}
If $c=c_*$, then \eqref{s3t6e1-1} gives
\[
\widetilde M_n(z;b,c;q)
=(z-bq^2)\widetilde M_{n-1}(z;bq^2,c/q^2;q).
\]
The two polynomials in the first assertion then coincide, so their
zeros do not strictly interlace. Moreover, $F_n^q(b,c)=bq^2=z_{1,n}$.\end{remark}

Finally, we compare the zeros of the degree $n$ quasi-orthogonal polynomial $\widetilde{M}_{n}\left(z; b, c; q\right)$ with the zeros of the degree $n-2$ orthogonal polynomial $\widetilde{M}_{n-2}\left(z; bq^2, c/q^2; q\right)$.

\begin{theorem}\label{s3t7}
    Let $n\ge3$, $0<q<1$, $q^{-1}<b<q^{-2}$, and $c>0$. Suppose that
$\widetilde M_n(z;b,c;q)$, which is quasi-orthogonal of order one with respect to the weight $w^q(x;bq,c/q)$, and
$\widetilde M_{n-2}(z;bq^2,c/q^2;q),$ which is orthogonal with respect to the weight $w^q(x;bq^2,c/q^2)$, have no common zeros.
Let $G_n^q(b,c)=\displaystyle\frac{q+c(1-bq)}{q^n}$.
Then the zeros of
$\widetilde M_n(z;b,c;q)$
and $(z-G_n^q(b,c))\widetilde M_{n-2}(z;bq^2,c/q^2;q)$
that lie in $(1,\infty)$ strictly interlace.
Moreover,
\begin{enumerate}
  \item[(i)] If $G_n^q(b,c)<z_{1,n}$, then the zeros of
$\widetilde M_n(z;b,c;q)$
and $(z-G_n^q(b,c))(z-1)\widetilde M_{n-2}(z;bq^2,c/q^2;q)$
strictly interlace.
  \item[(ii)] If $z_{1,n}<G_n^q(b,c)<z_{2,n}$, then the zeros of
$\widetilde M_n(z;b,c;q)$
and $(z-G_n^q(b,c))\widetilde M_{n-2}(z;bq^2,c/q^2;q)$
strictly interlace.
\end{enumerate}
\end{theorem}

\begin{proof}
Write $K=\displaystyle\frac{c(1-q^{n-1})}{(1-bq^n)(c+q^n)}>0.$

To prove the required interlacing result, we employ the following $q$-contiguous relation,
which is derived using \eqref{s1e2}, \eqref{s3t1e1}, and \eqref{s3l1e1}:
\begin{multline}\label{s3t7e1a}
  \frac{c q^{n-1} (1-b q)+q^n \left(1-b q^n\right)}{\left(1-b q^n\right) \left(c+q^n\right)} \widetilde{M}_{n}\left(z;b, c; q\right)
  = \left(z-G_n^q(b,c)\right) \widetilde{M}_{n-1}\left(z;b, c; q\right) \\
  -K  \left(z-bq\right) \left(z-bq^2\right) \widetilde{M}_{n-2}\left(z;bq^2, c/q^2; q\right).
\end{multline}

For $i=2,\ldots,n$, we have $z_{i,n}>bq>1>bq^2$. Evaluating \eqref{s3t7e1a} at $z_{i,n}$ gives
\[
(z_{i,n}-G_n^q(b,c))\widetilde M_{n-1}(z_{i,n};b,c;q)
=K(z_{i,n}-bq)(z_{i,n}-bq^2)\widetilde M_{n-2}(z_{i,n};bq^2,c/q^2;q).
\]
Since $\widetilde M_{n}(z;b,c;q)$ and $\widetilde M_{n-2}(z;bq^2,c/q^2;q)$ have no common zeros, this identity implies $G_n^q(b,c)\ne z_{i,n}$.

Define the monic polynomial
\[
S_{n-1}(z)=(z-G_n^q(b,c))\widetilde M_{n-2}(z;bq^2,c/q^2;q).
\]
Then
\[
S_{n-1}(z_{i,n})
=
\frac{(z_{i,n}-G_n^q(b,c))^2}
{K(z_{i,n}-bq)(z_{i,n}-bq^2)}
\,\widetilde M_{n-1}(z_{i,n};b,c;q),
\qquad i=2,\ldots,n.
\]
The prefactor is strictly positive. Thus $S_{n-1}$ and $\widetilde M_{n-1}(z;bq,c/q;q)$ have the same sign at each of these zeros of $P_n$.

By the ordering \eqref{s3t2e3}, these signs alternate, so $S_{n-1}$ has a zero in each interval
\[
(z_{i,n},z_{i+1,n}),\qquad i=2,\ldots,n-1.
\]
Furthermore, \eqref{s3t2e3} shows that exactly $n-2$ zeros of $\widetilde M_{n-1}(z;bq,c/q;q)$ exceed $z_{2,n}$. Hence
\[
\operatorname{sgn}S_{n-1}(z_{2,n})
=\operatorname{sgn}\widetilde M_{n-1}(z_{2,n};b,c;q)
=(-1)^{n-2}.
\]
Since $S_{n-1}$ is monic of degree $n-1$, its sign sufficiently far to the left is $(-1)^{n-1}$. Therefore, it has an additional zero in $(-\infty,z_{2,n})$.

These $n-1$ disjoint intervals account for all its zeros, which are consequently simple. Denoting them by $t_1<\cdots<t_{n-1}$, we obtain
\begin{equation}\label{ordering}
t_1<z_{2,n}<t_2<z_{3,n}<\cdots<t_{n-1}<z_{n,n}.
\end{equation}
Since $z_{1,n}<1$, restricting this ordering to $(1,\infty)$ proves the first assertion.

For either additional assertion, $G_n^q(b,c)<z_{2,n}$. As $G_n^q(b,c)$ is a zero of $S_{n-1}$, the ordering \eqref{ordering} implies $G_n^q(b,c)=t_1$.
Thus, the zeros of $\widetilde M_{n}(z;b,c;q)$ and $\widetilde M_{n-2}(z;bq^2,c/q^2;q)$ interlace as
\[
z_{i,n}<y_{i-1,n-2}<z_{i+1,n},
\qquad i=2,\ldots,n-1.
\]
If $G_n^q(b,c)<z_{1,n}$, the full ordering is
\[
G_n^q(b,c)<z_{1,n}<1<z_{2,n}<y_{1,n-2}<z_{3,n}
<\cdots<y_{n-2,n-2}<z_{n,n},
\]
which proves (i).

If $z_{1,n}<G_n^q(b,c)<z_{2,n}$, the ordering becomes
\[
z_{1,n}<G_n^q(b,c)<z_{2,n}<y_{1,n-2}<z_{3,n}
<\cdots<y_{n-2,n-2}<z_{n,n},
\]
which proves (ii). \end{proof}

\section*{Acknowledgements}
The research of the first author (KJ) was partially funded by DSTI-NRF Centre of Excellence in Mathematical and Statistical Sciences (CoE-MaSS), South Africa.
 Opinions expressed and conclusions arrived at are those of the authors and are not necessarily to be attributed to the CoE-MaSS.

\bibliography{PKKJ1}
\end{document}